\documentclass[11pt]{article}
\usepackage{amsmath,amsfonts}
\usepackage{xcolor}
\usepackage{enumerate}
\usepackage[letterpaper,centering]{geometry}

\numberwithin{equation}{section}
\newtheorem{theorem}{Theorem}[section]
\newtheorem{lemma}[theorem]{Lemma}

\newtheorem{corollary}[theorem]{Corollary}
\newtheorem{definition}[theorem]{Definition}
\newtheorem{remark}[theorem]{Remark}

\newcommand{\footremember}[2]{%
    \footnote{#2}
    \newcounter{#1}
    \setcounter{#1}{\value{footnote}}%
}
\newcommand{\qed}{\nolinebreak\hfill\vbox{\hrule\hbox{\vrule\kern3pt\vbox{\kern6pt}\kern3pt\vrule}\hrule}}
\newenvironment{pf}{{\noindent\bf Proof.}}{\qed\newline}

\newcommand{\supp}{\operatornamewithlimits{supp}}
\newcommand{\R}{\mathbb{R}}
\allowdisplaybreaks[4]

\begin{document}

\title{The obstacle problem for a higher order fractional Laplacian}

\author{Donatella Danielli\footremember{DD}{School of Mathematical and Statistical Sciences, Arizona State University, Tempe, Arizona, USA. Email: ddanielli@asu.edu},  Alaa Haj Ali \footremember{AH}{School of Mathematical and Statistical Sciences, Arizona State University, Tempe, Arizona, USA. Email:alaa.haj.ali@asu.edu}
  \& Arshak Petrosyan\footremember{AP}{Department of Mathematics, Purdue University, West Lafayette, Indiana, USA. \newline Email: arshak@purdue.edu}
}
\date{}
\maketitle
\begin{abstract}
In this paper, we consider the obstacle problem for the fractional Laplace operator $(-\Delta)^s$  in the Euclidian space $\mathbb{R}^n$ in the case where $1<s<2$. As first observed in \cite{Y}, the problem can be extended to the upper half-space $\mathbb{R}_+^{n+1}$ to obtain a thin obstacle problem for the weighted biLaplace operator $\Delta^2_b U$, where $\Delta_b U=y^{-b}\nabla \cdot (y^b \nabla U)$. Such a problem arises in connection with unilateral phenomena for elastic, homogenous, and isotropic flat plates.  We establish the well-posedness and $C_{loc}^{1,1}(\R^n) \cap H^{1+s}(\R^n)$-regularity of the solution. By writing the solutions in terms of Riesz potentials of  suitable local measures, we can base our proofs on tools from potential theory, such as a continuity principle and a maximum principle. Finally, we deduce the regularity of the extension problem to the higher dimensional upper half space. This gives an extension of Schild's work in \cite{Sc1} and \cite{Sc2} from the case $b=0$ to the general case $-1<b<1$.
\end{abstract}

\noindent {\bf Mathematics Subject Classification:} 35B65, 35R35, 35J35, 35R11\\

\noindent {\bf Keywords: } Free boundary problems, Variational Methods, Biharmonic operator, Higher order fractional Laplacian.\\

\noindent{\bf Statements and Declaration: }The authors have no relevant financial or non-financial interests to disclose.\\

\newpage

\tableofcontents

\section{Introduction}
Obstacle-type problems for second order uniformly elliptic operators have attracted  great interest over the years, and they continue to be explored, due the wide range of applications in elasticity, biology, and engineering,  just to name a few. In the past two decades there has been a resurgence of interest in obstacle problems associated with a class of nonlocal operators, whose prototype is the fractional Laplacian $(-\Delta)^s$ with order $s\in (0,1)$, driven by their relevance in financial mathematics. On the other hand, much less is known on obstacle-type problems associated with the bi-harmonic operator, or with higher order of the fractional Laplacian.

In this work we are concerned with the regularity of solutions to an obstacle problem associated with the fractional Laplacian $(-\Delta )^s u$ for $s\in (1,2)$, in the Euclidean space $\R^n$, with $n \geq 1$. In its simplest formulation, given a sufficiently regular function $\psi:\R^n\to\R$, known as the \emph{obstacle}, we want to find a function $u$ that satisfies
\begin{equation}\label{prb1}
(-\Delta)^s u\geq 0,\quad u\geq \psi,\quad (-\Delta)^s u (u-\psi)=0\quad\text{in }\R^n,
\end{equation}
with zero boundary condition at infinity
$$
u(x),\ \nabla u(x)\to 0\quad\text{as }|x|\to\infty.
$$
These conditions have to be understood in appropriate weak sense. The unique solution is obtained by minimizing the energy functional
\begin{equation}\label{prb1_fnct}
I_0[v]:=\int_{\R^n} \int_{\R^n} \frac{|\nabla v(x)-\nabla v(z)|^2}{|x-z|^{n-2+2s}} dx dz,
\end{equation}
(which is the square of the seminorm $[v]_{\dot H^{s}}$, see Section~\ref{prelim} for more details) over the space of functions
\begin{equation}\label{prb1_admissible}
\mathcal{B}_0:=\{v \in \dot{H}^{s}(\R^n)\mid v\geq \psi\},
\end{equation}
where $\dot{H}^s(\R^n)$ is the homogeneous fractional Sobolev space, which is the closure of $C^\infty_0(\R^n)$ with respect to $[\cdot]_{\dot  H^s}$.

We also consider another version of the problem, where besides the obstacle $\psi$, we are given a bounded domain $\Omega\in\R^n$ with sufficiently regular boundary. We want to find a function $u\in \dot{H}^s(\R^n)$ such that
\begin{equation}\label{prb2_original}
(-\Delta)^s u\geq 0,\quad u\geq \psi,\quad (-\Delta)^s u \ (u-\psi)=0\quad\text{in }\Omega,
\end{equation}
together with ``boundary'' conditions
$$
u=g\quad\text{on }\R^n\setminus\Omega,
$$
where $g\in \dot{H}^s(\R^n)$ satisfies the compatibility condition $g\geq\psi$ on $\R^n\setminus\Omega$. If we additionally assume that $g\in C^{2s+\epsilon}_{\rm loc}\cap \dot H^{2s}$, then replacing $u$ with $u-g$ and $\psi$ with $\psi-g$, and letting $f=(-\Delta)^sg\in L^2(\R^n)\cap L^\infty_{\rm loc}(\R^n)$, we can equivalently state the above problem as
\begin{equation}\label{prb2}
(-\Delta)^s u\geq f(x),\quad u\geq \psi,\quad ((-\Delta)^su-f) (u-\psi)=0\quad\text{in }\Omega,
\end{equation}
with boundary conditions
$$
u=0\quad\text{on }\R^n\setminus\Omega.
$$
The solution of (\ref{prb2}) can be obtained by minimizing the functional
\begin{equation}\label{prb2_fnct}
I[v]:=\frac{1}{2k(n,s)}\int_{\R^n} \int_{\R^n} \frac{|\nabla v(x)-\nabla v(z)|^2}{|x-z|^{n-2+2s}} dx dz -\int_{\R^n} f(x) v(x) dx
\end{equation}
over the set
\begin{equation}\label{prb2_admissible}
\mathcal{B}:=\{v \in \dot{H}^{s}(\R^n)\mid v=0\text { on } \R^n \setminus\Omega\text { and } v\geq \psi\}.
\end{equation}
The constant $k(n,s)$ comes from the equivalency (\ref{inner_product}) below.

Inspired by Schild's work in \cite{Sc1} and \cite{Sc2}, we study the
regularity of the minimizers of $I_0[\cdot]$ and $I [\cdot]$ using
techniques from potential theory.  In  \cite{Sc2} (resp. \cite{Sc1})
Schild studies the thin obstacle problem for the poly-harmonic
operator $\Delta^{m} u=0$, for integers $m\geq 2$,  in a bounded
domain $\Omega$ of the Euclidean space for the two dimensional case
(resp. for dimension greater than two). Among other results, the
author establishes the $C_{\rm loc}^{1,1} \cap H_{\rm loc}^{3}$-regularity of the solution.

To apply regularity results from potential theory, given a ball $B'_{\rho}$ in $\R^n$ and a measure $\mu_{\rho}$ on $\mathbb{R}^n$  supported on $B'_{2\rho}$, we aim to represent the minimizer and its second order partial derivatives on $B_{\rho}'$ in terms of the Riesz potential of $\mu_{\rho}$ of order $s$
$$\mathcal{R}_{\mu_{\rho}}:=\phi_{s}*\mu_{\rho},$$
where $\phi_{s}$ is the fundamental solution of the fractional Laplace operator $(-\Delta)^s$ as defined in (\ref{fundamental_fractional}).
As a starting point, for problem (\ref{prb1}) we show that the minimizer of $I_0[\cdot]$ is equal to $\mathcal{R}_{\mu_{\rho}}$ almost everywhere on $\R^n$. Therefore, for problem (\ref{prb1}), we are able to derive the desired representation on the whole space $\R^n$. However, this is not the case with problem (\ref{prb2}), where boundary data are imposed on $\R^n\setminus\Omega$. Instead, for problem (\ref{prb2}) we show that $u= \mathcal{R}_{\mu_{\rho}}+R_{\rho}$ on $B'_{\rho}$, where $R_{\rho}$ is a weak solution of $(-\Delta)^s R_{\rho}=f$ on $B'_{\rho}$. In order to obtain an estimate for $\mathcal{R}_{\mu_{\rho}}$ in terms of the Sobolev norms of the minimizer, we need to localize our problem by applying an extension result from \cite{Y} (see also \cite{CM22}), in the spirit of the well-known Caffarelli-Silvestre \cite {CS} extension for the case $0<s<1$.
In  \cite{CS} the authors characterize the fractional Laplacian of a function $u$ on $\R^n$ as the Dirichlet-to-Neumann map for a function $U$ satisfying, in some weak sense, the degenerate elliptic equation
$$
L_a U:=y^{-a}\nabla \cdot (y^a \nabla U)=0
$$
in the upper half space $\R^n\times\R_+$ (here $a=1-2s$). More precisely, they show that
$$(-\Delta)^{s} u=-C_{n,a} \lim_{y \to 0^{+}} y^a U_y(x,y)$$ for a suitable positive constant $C_{n,a}$. With this characterization, one can avoid working directly with the fractional Laplacian to understand the obstacle problem associated with $(-\Delta)^s$, $0<s<1$, and study instead a local extension problem, known as the \textit{thin obstacle problem}

\begin{equation}\label{Signorini}
\left\{
\begin{array}{ll}
L_{a} U =0 & \text { for } (x,y) \in B_1^{+}\\
U(x,y)=h(x,y) & \text { for } (x,y) \in (\partial B_1)^+ \\
U(x,0) \geq \psi(x) & \text { in } \bar{B_1'}\\
\lim_{y \to 0^+} y^a U_{y}(x,y) \leq 0 & \text { in }  \bar{B_1'}\\
\big(U(x,0)-\psi(x)\big) \lim_{y \to 0} y^a U_{y}(x,y)=0& \text { in }  \bar{B_1'}.\\
\end{array}
\right.
\end{equation}
Here  $B_1$ is the unit ball in $\R^{n+1}$, $B_1^+=B_1 \cap \{y>0\}$, $(\partial B_1)^+=\partial B_1 \cap \{y>0\}$, and $h$ represents some Dirichlet boundary data. The special case $s={1}/{2}$ is referred to as the \textit{Signorini problem}.
The thin obstacle problem has been extensively studied. Just to mention a few salient results, in \cite{C}, Caffarelli proves that a solution of the Signorini problem is $C^{1,\alpha}$ up to $y=0$ for some $0<\alpha<{1}/{2}$. This result was improved upon by Athanasopoulos and Caffarelli in \cite{AC}, where they achieved the optimal $C^{1,{1}/{2}}_{\rm loc}$-regularity of the solution. In the case of vanishing obstacle, Athanasopoulos, Caffarelli and Salsa in \cite {ACS} analyzed the regularity of the free boundary near the so-called regular free boundary points. In 2008, Caffarelli, Salsa and Silvestre in \cite{CSS} extended the optimal regularity results and the analysis of the regular part of the free boundary to all $0<s<1$ and non-vanishing obstacles.

The Caffarelli-Silvestre characterization of the fractional Laplacian as a Dirichlet-to-Neumann map was generalized in \cite{Y} to positive, non-integer orders of the fractional Laplace operator.  More specifically, for $1<s<2$, let  $b=3-2s$ and let $\Delta_b$ be the operator defined by
\begin{equation}\label{b-Laplace}
\Delta_b U=y^{-b} \nabla \cdot (y^b \nabla U).
\end{equation}
For $\Omega$ an open domain in $\R^{n+1}$, let $H^2 (\Omega,y^b)$ be the weighted Sobolev space equipped with the norm
 \begin{equation}\label{weight_norm}
 \|U\|^2_{H^2 (\Omega,y^b)}=\|y^{\frac{b}{2}} \Delta_b U\|^2_{L^2 (\Omega)}+\|y^{\frac{b}{2}} \nabla U\|^2_{L^2 (\Omega)}+\|y^{\frac{b}{2}}  U\|^2_{L^2 (\Omega)}.
 \end{equation}
We recall the following definition of a weak solution of the $b$-biharmonic equation $\Delta_b^2U=0$.
 \begin{definition}\label{bi-harmonic*}
 We say that $U \in H^{2} (\Omega,y^b) $ is a weak solution of the b-biharmonic equation $\Delta_b^2U=0$ in $\Omega$ if
 \begin{equation}\label{bi-harmonic}
 \int_{\Omega} y^{b}\ \Delta_b U \Delta_b \phi =0 \text { for all functions }\phi \in H_0^{2}(\Omega,y^b).
 \end{equation}
\end{definition}
The following result holds.
\begin{theorem}\textup{\cite[Theorem 3.1]{Y}}\label{extension*}
Let $1<s<2$ be a non-integer and let $b=3-2s$.
Suppose $U$ is a function in $H^{2} (\R_{+}^{n+1},y^b)$ satisfying
\begin{equation*}
\left\{
\begin{array}{ll}
 \Delta^2_b U(x,y)=0 & \text { in } \R^n \times \R_{+}\\
U(x,0)=u(x) & \text { for all } x \in \R^n\\
\lim_{y \to 0} y^b U_{y}(x,y)=0 & \text { for all } x \in \R^n\\
\end{array}
\right.
\end{equation*}
for some $u\in \dot{H}^{s}(\R^n)$.
Here  $\Delta_b$ is as  in (\ref{b-Laplace}) and the Dirichlet boundary condition is satisfied in the sense of traces.
Then,
\begin{equation}\label{D_to_N}
(-\Delta)^s u(x)=C_{n,b} \lim_{y \to 0} y^b\frac{\partial}{\partial y} \Delta_b U(x,y).
\end{equation}
\end{theorem}

Applying Theorem \ref{extension*}, we see that problem (\ref{prb2}) is equivalent to a boundary obstacle problem in $ B_1^{+}$, with an obstacle living on $B_1'$.
For the limiting case $b=0$, the operator $\Delta_b^2$ becomes the biharmonic operator - this problem is the object of Schild's work. Therefore, our paper can also be seen as a generalization of such work to the case of the operator $\Delta^2_b$ for $ b \in (0,1)$.

Our paper is organized as follows. In Section \ref{prelim}, we define the fractional Laplacian $(-\Delta)^s, \ \  (1<s<2)\ $, and we  introduce the concept of weak solution in $\dot{H}^s(\R^n)$ of equations involving such an operator.  In Section \ref{existence_uniqueness}, we establish the existence and uniqueness of a minimizer $u_0$ (resp. $u$) of the functional $I_0[\cdot]$ (resp. $I[\cdot]$) over the set $\mathcal{B}_0$ (resp. $\mathcal{B}$). Next, we derive the variational inequality satisfied by $u_0$ and $u$, and we introduce the corresponding measures $\mu_0$ and $\mu$, respectively. In Section \ref{problem1}, we show that the minimizer $u_0$ of the functional $I_0[\cdot]$ is equal a.e. on $\R^n$ to some Riesz potentials of the measure $\mu_0$. In section \ref{extension_problem}, we extend a function $u$ on $\R^n$ to a b-biharmonic function $U$ on the upper half (n+1)-dimensional space  by means of a convolution with the Poisson kernel which we derive. In Section \ref{problem2}, we aim to obtain a local a.e. representation of the minimizer $u$ of $I[\cdot]$ as a sum of Riesz potential and a function $R_{\rho}$ satisfying $(-\Delta)^s R_{\rho}=f$. To estimate the function $R_\rho$ in terms of the solution $u$, we need to extend $u$ to the higher dimensional space in the spirit of the extension Theorem \ref{extension*}. The results of the Section 7 and 8 apply to both of the minimizers $u_0$ and $u$, and they aim at establishing their regularity. In Section \ref{C1alpha_regularity},  with the representations of Sections \ref{problem1} and \ref{problem2} in hand, we use techniques from potential theory to establish the $C^{1,1}_{\rm loc}$-regularity of the minimizers. More specifically, using a continuity principle from potential theory, we first prove the continuity of the minimizers. Then, with an application of Green's formula and an adaptation of a maximum principle from potential theory, we are able to establish the $C_{\rm loc}^{1,1}$-regularity. In Section \ref{H1+s_regularity}, using a result which relates the $L^2$-norm of the gradient of Riesz potential and the Riesz energy of a measure, we establish the $H_{\rm loc}^{1+s}$-regularity of a minimizer. Finally, in Section 9 we deduce the regularity of the extension function $U$ of $u$ to the higher dimensional upper half space $\R_+^{n+1}$, hence extending Schild's work from the case $b=0$ to the general case $-1<b<1$.

Throughout the paper, for a given function $W$ defined on $\R^d$ for some dimension $d$, $\nabla W$ (resp. $\Delta W$) will denote the gradient of $W$ (resp. the Laplacian of $W$) on $\R^d$. Moreover, for a function $W(x,y)$ defined on $\R^{n} \times \{y \geq 0\}$, $\nabla_x W$ (resp. $\Delta_x W$) will denote the gradient of $W$ (resp. the Laplacian of $W$) with respect to the variable $x$.

\section{\boldmath Fractional Sobolev spaces and the fractional Laplacian of order $1<s<2$}\label{prelim}

For $1<s<2$, the homogeneous fractional Sobolev space $\dot H^s(\R^n)$ is defined as the closure of $C_0^\infty(\R^n)$ with respect to the seminorm
\begin{align}\label{semi-norm}
[v]_{\dot H^s}&:=\left(\int_{\R^n} \int_{\R^n} \frac{|\nabla v(x)-\nabla v(z)|^2}{|x-z|^{n-2+2s}} dx dz\right)^{1/2}\\
&=\left(k(n,s)\int_{\R^n}|\xi|^{2s}|\hat v(\xi)|^2d\xi\right)^{1/2},
\end{align}
where $\hat v(\xi)$ is the Fourier transform of $v$ and $k(n,s)$ is a positive constant depending on the dimension $n$ and the exponent $s$ only.
The space $\dot{H}^{s}(\R^n)$ is a Hilbert space when endowed with the inner product
\begin{align}
\label{inner_product}
\langle f,g\rangle_{\dot{H}^s}&=\int_{\R^n} \int_{\R^n} \frac{(\nabla f(x)-\nabla f(z)) \cdot (\nabla g(x)-\nabla g(z))}{|x-z|^{n-2+2s}} dx dz\\
&=k(n,s)\int_{\R^n}|\xi|^{2s} \hat f(\xi)\overline{\hat g(\xi)}d\xi
\end{align}

For a function $v \in \dot{H}^s(\R^n)$ and $1<s<2$ we define the \emph{fractional Laplacian} $(-\Delta)^sv$ as a tempered distribution on $\R^n$ such that
\begin{equation}\label{frac}
\widehat{(-\Delta)^s v}(\xi)=|\xi|^{2s}\hat v(\xi),\quad \xi\in\R^n
\end{equation}
in terms of Fourier transforms. Equivalently, it can be given in the weak form as
\begin{equation}\label{frac_weak}
\langle (-\Delta)^s v,\phi \rangle=\frac{1}{k(n,s)}\langle v,\phi\rangle_{\dot H^s(\R^n)},\quad\text{for any }\phi\in C^\infty_0(\R^n).
\end{equation}


In our proofs, we will also use the following Sobolev embedding theorems. If $u\in \dot H^s(\R^n)$, $1<s<2$, then $u\in H^1_{\rm loc}(\R^n)$ and moreover $u\in L^{\frac{2n}{n-2s}}(\R^n)$, $\nabla u\in L^{\frac{2n}{n+2-2s}}(\R^n)$ and
\begin{equation}\label{Sobolev_ineq}
\|u\|_{L^{\frac{2n}{n-2s}}(\R^n)}+\|\nabla u\|_{L^{\frac{2n}{n+2-2s}}(\R^n)}\leq C_{n,s}[u]_{\dot H^s}.
\end{equation}
For a domain $\Omega \subset \subset \R^n$, we say that $u \in \dot{H}^s(\Omega)$ if its extension $$\tilde{u}:=\left\{\begin{array}{ll} u & x \in \Omega\\ 0 & x \in \R^n \setminus \Omega
\end{array}\right.$$
is in $\dot{H}^s(\R^n)$.
\section{Existence and Uniqueness of solutions}\label{existence_uniqueness}
In this section  we establish the existence and uniqueness of a minimizer $u_0$ (resp. $u$) of the functional $I_0[\cdot]$ (resp. $I[\cdot]$) over the set $\mathcal{B}_0$ (resp. $\mathcal{B}$). In addition, we derive the variational inequalities satisfied by $u_0$ and $u$.

\begin{theorem}[Existence and Uniqueness of  Minimizers]\label{Existence}
There exists a unique minimizer $u_0$ (resp. $u$) of the functional $I_0[\cdot]$ (resp. $I[\cdot]$) over the admissible set $\mathcal{B}_0$ (resp. $\mathcal{B}$).
\end{theorem}
\begin{pf}
Concerning $u_0$, first we note that the admissible set $\mathcal{B}_0$ is non-empty since $\psi \in \mathcal{B}_0$. Now let $\{u_m\}$ be a minimizing sequence for $I_0[\cdot]$. Then, $\{u_m\}$ is bounded in the $H^s$-norm, and therefore it has a  subsequence which converges weakly in $\dot{H}^s(\R^n)$ to a function $u_0\in\mathcal{B}_0$.
By the property of weak convergence, we immediately conclude that $u_0$ is a minimizer of $I_0[\cdot]$ over the set $\mathcal{B}_0$.

An analogous argument shows that the minimizing sequence of $I[\cdot]$ (still denoted by $\{u_m\}$), has a  subsequence which converges weakly in $\dot{H}^s(\R^n)$ to a function $u\in\mathcal{B}$. The Sobolev inequality (\ref{Sobolev_ineq}) implies that $\{u_m\}$ is bounded in $ H^{1}(\R^n)$. Moreover, since $u_m \equiv 0$ in $\R^n \setminus B'_1$, we apply the Rellich-Kondrachov compactness theorem for $H_0^{1}(B'_1)$ to conclude that $\{u_m\}$ converges strongly in $L^2(B'_1)$. Passing to the weak limit in $I[\cdot]$ we conclude that $u$ is indeed a minimizer of $I[\cdot]$. The uniqueness of both of minimizers follows from the convexity of the functionals $I_0[\cdot]$ and $I[\cdot]$.
\end{pf}

\noindent
In the next two lemmas, we obtain the variational inequalities satisfied by the minimizers $u_0$ and $u$.

\begin{lemma}\label{prb1_varineq*}
The minimizer $u_0$ of the functional $I_0[\cdot]$ over the set $\mathcal{B}_0$ satisfies the variational inequality
\begin{equation}\label{prb1_varineq}
\langle (-\Delta)^s u_0, v-u_0\rangle=\frac1{k(n,s)}\langle u_0, v-u_0\rangle_{\dot H^s}\geq 0,\quad\text{for all }v \in \mathcal{B}_0.
\end{equation}
In particular,
\begin{equation*}
\mu^0:=(-\Delta )^s u_0 \geq 0
\end{equation*}
is a nonnegative Radon measure on $\R^n$.
\end{lemma}

\begin{pf}
For $v \in \mathcal{B}_0$, write $v=u_0+\phi$. Then $u_0+\epsilon \phi$ is an admissible function for all $\epsilon>0$.
By the minimality of $u_0$ we have
\begin{equation*}
0\leq \langle u_0+\epsilon \phi, u_0+\epsilon
\phi\rangle_{\dot{H}^{s}(\R^n)} -\langle
u_0,u_0\rangle_{\dot{H}^{s}(\R^n)}=2 \epsilon \langle u_0,\phi\rangle_{\dot{H}^{s}(\R^n)} + \epsilon^2 \langle\phi,\phi\rangle_{\dot{H}^{s}(\R^n)}.
\end{equation*}
Dividing by $\epsilon$ and then letting $\epsilon \to 0$ we obtain
$$\langle u_0,\phi\rangle_{\dot{H}^{s}(\R^n)} \geq 0.$$
\end{pf}

\noindent
We can derive the following variational inequality for the minimizer $u$ of $I[\cdot]$ in a similar way .

\begin{lemma}\label{prb2_varineq*}
The minimizer $u$ of the functional $I[\cdot]$ over the set $\mathcal{B}$ satisfies the variational inequality
$$
\langle (-\Delta)^s u-f, v-u\rangle=\frac1{k({n,s)}}\langle u, v-u\rangle_{\dot H^s}-\langle f, v-u\rangle\geq 0,\quad\text{for all }v \in \mathcal{B}.
$$
In particular,
\begin{equation*}
\mu:=(-\Delta )^s u -f
\end{equation*}
is a nonnegative Radon measure on $\Omega$.
\end{lemma}

\section{The problem without boundary values - Problem (\ref{prb1})}\label{problem1}

In the remainder of the paper, $\phi_s$ represents the fundamental solution of the fractional Laplacian $(-\Delta)^s$, which is given by
 \begin{equation}\label{fundamental_fractional}
 \left\{\begin{array}{lll}
\phi_s(x)=\frac{1}{|x|^{n-2s}},&\quad x \in \R^n \setminus \{0\}   & \text { when $n \neq 1,2,3$ or $s \neq \frac{3}{2}$},\\
\phi_{\frac{3}{2}}(x)=-\log|x|,&\quad x \in \R^n \setminus \{0\}    & \text { when $n=3$,}\\
\phi_{\frac{3}{2}}(x)=-|x|,&\quad x \in \R^n \setminus \{0\}    & \text { when $n=2$,}\\
\phi_{\frac{3}{2}}(x)=|x|^2\left(\log|x|-\frac{3+n}{2n+2}\right),&\quad x \in \R^n \setminus \{0\}    & \text { when $n=1$.}
\end{array}\right.
\end{equation}
To study the regularity of the minimizers $u_0$ and $u$, we use techniques from potential theory. We start with representing the minimizers in terms of Riesz potentials of measures. For problem \eqref{prb1}, we are able to obtain an a.e.  representation of $u_0$ as the Riesz potential of order $2s$ of the measure $\mu^0$.

\begin{lemma}\label{u0_rep*}
The identity
\begin{equation}\label{u0_rep}
u_0=\phi_s*\mu^0
\end{equation}
 holds true for a.e.\ $x$ in $\R^n$.
\end{lemma}

\begin{pf}
First of all, it is not hard to see that $\phi_s*\mu^0$ and $
\Delta \phi_s*\mu^0$ are in $L^1_{loc}(\mathbb{R}^n)$, and therefore finite almost everywhere. Now, let $\zeta \in C_0^{\infty}(\R^n)$.
We have
\begin{equation}\label{est1}
\int_{\R^n} u_0(x)\zeta(x) dx=\int_{\R^n} \hat{u}_0(\xi)\overline{\hat{\zeta}(\xi)}d\xi=\int_{\R^n}|\xi|^{2s}\hat{u}_0(\xi)\overline{\hat{\phi}_s(\xi)\hat\zeta(\xi)}=\int_{\R^n} (\phi_s*\zeta)(x) d\mu^0(x).
\end{equation}
 Moreover, applying Fubini's theorem we have
 \begin{equation}\label{est2}
 \begin{split}
 \int_{\R^n} (\phi_s*\zeta)(x) d\mu^0(x)&=\int_{\R^n} \left(\int_{\R^n} \phi_s(y-x)\zeta(x) dx \right)d\mu^0(y)\\&=\int_{\R^n} \left(\int_{\R^n} \phi_s(y-x) d\mu^0(y)\right)\zeta(x) dx
 =\int_{\R^n}  (\phi_s*\mu^0)(x)\,\zeta(x)dx.\\
 \end{split}
 \end{equation}
Combining (\ref{est1}) with (\ref{est2}) yields
\begin{equation*}
 u_0=\phi_s*\mu^0\quad\text {a.e.\ in } \R^n.
 \end{equation*}
\end{pf}
\section{The extension problem}\label{extension_problem}
To study the regularity of problem (\ref{prb2}), we extend the minimizer $u$ to the upper half-space with a function $U\in H^2(\R^n,y^b)$  satisfying
\begin{equation}\label{extension_prb}
\left\{\begin{array}{ll}
\Delta_b^2 U=0 & \text { for } (x,y) \in \R^{n+1}_+,\\
\lim_{y \to 0^+} y^b U_y (x,y)=0 &\text { for } (x,y) \in \R^{n+1}_+,\\
U(x,0)=u(x)  & \text { for }x \in \R^n.\\
\end{array}\right.
\end{equation}
HJere $b=3-2s$, and the PDE  is satisfied in the weak sense (\ref{bi-harmonic}). We say that $U$ is the \emph{b-biharmonic extension} of $u$ to the upper half-space $\R_+^{n+1}$.  We begin with deriving the fundamental solution of the b-biLaplace equation $\Delta_b^2U=0$ at the origin.

\subsection{Fundamental solution at the origin}
For $s\in (1,2)$ let $b=3-2s\in (-1,1)$, or equivalently $s=\frac{3-b}{2}$. First, we make the following observation. Suppose $w:\R^n \times \R^{1+b} \to \R$ is radially symmetric in the  second variable $y$, e.g. $w(x,y_1)=w(x,y_2)$ whenever $|y_1|=|y_2|$. Denoting $r:=|y|$, a direct computation shows that
\begin{equation*}
\Delta^2 w (x,y)= \Delta_b^2 w(x,r).
\end{equation*}
Thus, following the line of thought in \cite{CS}, we can think of the problem (\ref{extension_prb}) as the biharmonic extension problem of $u$ in $1+b$ additional dimensions.

With this observation, we see that the fundamental solution at the origin of the biLaplace equation in $n+1+b$ dimensions  is indeed a fundamental solution for the b-biLaplace equation. From the Appendix in \cite{Sc1}, we know that the latter reads, for $ x \in \R^n \setminus \{0\}$, as

 \begin{equation}\label{fundamental}
 \left\{\begin{array}{lll}
E_s(x,y)=\frac{1}{\left(|x|^2+|y|^2\right)^{\frac{n-2s}{2}}},&\quad n \neq 1,2,3\text{ or } s \neq {3}/{2},\\
E_{{3}/{2}}(x,y)=-\log\left(|x|^2+|y|^2\right)^{\frac{1}{2}},&\quad n=3\\
E_{{3}/{2}}(x,y)=-\left(|x|^2+|y|^2\right)^{\frac{1}{2}},&\quad n=2,\\
E_{{3}/{2}}(x,y)=\left(|x|^2+|y|^2\right)\left(\log\left(|x|^2+|y|^2\right)^{\frac{1}{2}}-\frac{3+n}{2n+2}\right),&\quad n=1.\\
 \end{array}\right.
\end{equation}
Recalling (\ref{fundamental_fractional}), we notice that $E_{s}(x,0)$ is the fundamental solution of the fractional Laplacian $(-\Delta)^s$.

\subsection{Poisson Formula}
Let $b=3-2s$ as before. Now we derive the Poisson kernel for the b-biharmonic operator $\Delta_b^2$ on the upper half space, and therefore we obtain a Poisson representation formula for the extension problem (\ref{extension_prb}).

\begin{lemma}\label{l.Poisson_kernel}
The function
\begin{equation}\label{Poisson_kernel}
P_s(x,y)=-C_{n,s} \frac{y^{2s}}{ (|x|^2+|y|^2)^{\frac{n+2s}{2}}}
\end{equation}
is the Poisson kernel for the b-biLaplace equation  $\Delta_b^2 u=0$. That is, $\Delta_b^2 P_s(x,y)=0$ for all $(x,y)$ such that $y>0$, and $\lim_{y \to 0^+} P_s(x,y)=-\delta_0(x)$, where $\delta_0$ denotes  Dirac's delta. Moreover, $P_s(x,y)$ satisfies
\begin{equation}\label{P_to_E}
P_s(\cdot,y)=C_{n,s} (-\Delta)^{s} E_s (\cdot,y) \quad\text{for all $y>0$}.
\end{equation}
In particular, we have the following Poisson formula
\begin{equation}\label{poisson_representation}
U(x,y)=\int_{\R^n} P_s(x-\xi,y) u(\xi) d\xi=\left(E_{s}(\cdot,y)*(-\Delta)^s u\right)(x) \quad\text{for all $y>0$}.
\end{equation}
\end{lemma}
\begin{pf}
By direct computations we obtain
\begin{equation*}
\begin{split}
\Delta_b \left(\frac{y^m}{r^k}\right)&=(m^2+bm-m)\frac{y^{m-2}}{r^k}+k(k-n-2m+1-b)\frac{y^{m}}{r^{k+2}}\\
\end{split}
\end{equation*}
and
\begin{equation*}
\begin{split}
\Delta^2_b \left(\frac{y^m}{r^k}\right)
&=(m^2+bm-m)\left((m^2+bm-5m-2b+6)\frac{y^{m-4}}{r^{n+m}}+(n+m)(-2m+6-2b)\frac{y^{m-2}}{r^{n+m+2}}\right)\\
&+(n+m)(-m+1-b)(n+m+2)(-m+3-b)\frac{y^{m}}{r^{n+m+4}}.\\
\end{split}
\end{equation*}
Setting $k=m+n$ we get
\begin{equation*}
\begin{split}
\Delta^2_b \left(\frac{y^m}{r^k}\right)&=(m^2+bm-m)\left((m^2+bm-5m-2b+6)\frac{y^{m-4}}{r^{n+m}}+(n+m)(-2m+6-2b)\frac{y^{m-2}}{r^{n+m+2}}\right)\\
&+(n+m)(-m+1-b)(n+m+2)(-m+3-b)\frac{y^{m}}{r^{n+m+4}}.\\
\end{split}
\end{equation*}
Now we notice that for $m=3-b$, we have $\Delta^2_b \left(\frac{y^m}{r^k}\right)=0$. This shows that $P_s(x,y)$ as defined in (\ref{Poisson_kernel}) satisfies $\Delta_b^2 P_s =0$.
Moreover, we can easily check that $\lim_{y \to 0^+} P_s(x,y)=-\delta_0(x)$, where $\delta_0(x)$ is the Dirac delta function at 0. Therefore, $P_s(x,y)$ is indeed the Poisson kernel.

It remains to prove (\ref{P_to_E}). For this, we write
\begin{equation*}
E_s(\cdot,y)=\phi_s*P_s(\cdot,y).
\end{equation*}
Thus,
\begin{equation*}
\widehat{ E_s(.,y)}=\widehat{\phi_s}\widehat{P_s(\cdot,y)}=\frac{1}{|\xi|^{2s}} \widehat{P_s(\cdot,y)}.
\end{equation*}
This shows that
\begin{equation*}
\widehat{(-\Delta)^{s} E_{s}(.,y)}=|\xi|^{2s} \widehat{E_s(\cdot,y)}=\widehat{P_s(\cdot,y)} \quad\text{for all }y>0,
\end{equation*}
and yields (\ref{P_to_E}).
\end{pf}

\section{The problem with boundary values - Problem \textbf{(\ref{prb2})}}\label{problem2}

Throughout this section, $u$ will be the minimizer of the functional $I[\cdot]$. Using Theorem \ref{extension*} and the results of Section \ref{extension_problem}, we aim to obtain a local representation of $u$. Given the extension result of Theorem \ref{extension*}, we know that the extension function $U(x,y)$  of $u$ solves
\begin{equation}\label{prb2_extension}
\left\{
\begin{array}{ll}
\Delta^2_b U =0 & \text { for } (x,y) \in B_1^{+}\\
U=0 , \  U_{\nu}=0 & \text { for } (x,y) \in (\partial B_1)^+ \\
\lim_{y \to 0} y^b U_{y}(x,y)=0  & \text { for all } x \in B_1'\\
U(x,0) \geq \psi(x) & \text { in } B_1'\\
\lim_{ y \to 0^+} y^b \frac{\partial}{\partial y} \Delta_b U(x,y) \geq f(x)  & \text{ in } B'_1\\
\lim_{ y \to 0^+} y^b \frac{\partial}{\partial y} \Delta_b U(x,y) =f(x)  & \text{ in } B'_1 \cap \{u>\psi\},\\
\end{array}
\right.
\end{equation}
We recall that  $B_1$ is the unit ball in $\R^{n+1}$, $B_1^+=B_1 \cap \{y>0\}$, $(\partial B_1)^+=\partial B_1 \cap \{y>0\}$, and $\nu$ is the outer unit normal to the ball $B_1$. We will also denote $\overline{\R^{n+1}}:=\R^{n+1} \cap \{y \geq 0\}$
and $\overline{B^+}:=B \cap \{y \geq 0\}$ for a ball $B\subset \R^{n+1}$.
We will need the following approximation result, proven by Silvestre in \cite{S}.

\begin{lemma}\cite[Proposition 2.13 \& Proposition 2.15]{S}\label{lower_semi_cont}
For $0<\sigma<1$, let $\Gamma_\sigma$ be a $C^{1,1}$-function on $\mathbb{R}^n$ which coincides on $B'_1$ with the fundamental solution of the fractional Laplacian $(-\Delta)^\sigma$. Also, for $\epsilon >0$, let $\Gamma_{\sigma,\epsilon}=\frac{1}{\epsilon^{n-2\sigma}}\Gamma_\sigma\left(\frac{x}{\epsilon}\right)$ and
$\gamma_{\sigma,\epsilon}:=(-\Delta)^{\sigma}\Gamma_{\sigma,\epsilon}$. Then, the family of functions $\gamma_{\sigma,\epsilon}$ is an approximation of the identity, in the sense
\begin{equation}\label{a.e._approx_s}
h*\gamma_{\sigma,\epsilon} \to h\quad \text {\
  a.e. as}\quad\epsilon \to 0^+\quad\text {for all }h \in L^1.
\end{equation}
 Moreover,
in the case where the function $h$ satisfies $(-\Delta)^{\sigma} h \geq 0$,  then the convergence in (\ref{a.e._approx_s}) is monotone increasing, and therefore $h$ is equal a.e. to a lower semi-continuous function.
\end{lemma}

\noindent
Throughout the remainder of the paper, for $B'_{\rho}$ ball in $\R^n$ centered at the origin with radius $\rho$, we will denote by $\mu_{\rho}$ the restriction of the measure $\mu$ to $B'_{2\rho}$, that is \begin{equation}\label{measure_ball}
\mu_{\rho}:=\eta\mu,
\end{equation}
where $\eta \in C_0^{\infty}(B'_{2 \rho})$ is a cut-off function such that $\eta \equiv 1$ on $B'_{\rho}$ and $|\nabla \eta| \leq \frac{1}{\rho}$ on $B'_{2\rho}$.
We first obtain the following local representation of $u$ and $\Delta u$.

\begin{lemma}\label{rep_local*}
For a ball $B'_{\rho}$ in $\R^n$, there exists a weak solution $R_{\rho}(x)$ of $(-\Delta)^s R_{\rho}(x)= f$ in $B'_{\rho}$ such that
\begin{enumerate}[i)]
\item\begin{equation}\label{u_rep_local}
u(x)=\phi_{s}*\mu_{\rho}(x)+ R_{\rho}(x) \quad\text{for a.e.  $x \in B'_{\rho}$. }
\end{equation}
\item
\begin{equation}\label{lapxu_rep_local}
\Delta u(x)=\Delta \phi_{s}*\mu_{\rho}(x)+ \Delta R_{\rho}(x)\quad\text{for a.e. $x \in B'_{\rho}$. }
\end{equation}
\item $R_{\rho}$ satisfies the estimate
\begin{equation}\label{R_bnd}
\|R_{\rho}\|_{W^{2,\infty}(B'_{\rho})} \leq C(n,s,\rho) \left(\|y^{{b}/{2}} U\|_{L^2(B_{2 \rho})}  +\| y^{{b}/{2}} \nabla U\|_{L^2(B_{2 \rho})}+\|f\|_{L^{\infty}(B'_{2 \rho})} \right).
\end{equation}
\end{enumerate}
\end{lemma}
\begin{pf}
For $\rho>0$, define
\begin{equation}\label{R_rho_def}
R_{\rho}(x):=u(x)-\phi_{s}*\mu_{\rho}(x).
\end{equation}
\noindent
We aim to show that $R_{\rho}$ satisfies the equation $(-\Delta)^s R_{\rho}=f$ in $B'_{\rho}$.
First of all, it is not hard to see that $\phi_s*\mu_{\rho}$ and $
\Delta \phi_s*\mu_{\rho}$ are in $L^1_{loc}(\mathbb{R}^n)$, and therefore finite a.e.
For $\zeta \in C^{\infty}_0(B'_{\rho}) $, and $\eta$ the cut-off function supported in $B'(2\rho)$ defined right before (\ref{measure_ball}),
we compute
\begin{equation}\label{eq.interm}
\begin{split}
&\int_{\R^n}R_{\rho}(x) (-\Delta)^s \zeta(x) dx\\
 &=\int_{\R^n}u(x)(-\Delta)^s \zeta(x) dx -\int_{\R^n} (\phi_s*\mu_{\rho})(x)(-\Delta)^s \zeta(x) dx \\
&=\int_{\R^n}u(x)(-\Delta)^s \zeta(x) dx - \int_{\R^n} \left(\int_{\R^n}   \left((-\Delta)^s u(z)-f(z)\right) \eta(z) \phi_s(x-z) dz\right)  (-\Delta)^s \zeta(x) dx  \\
&=\int_{\R^n}u(x)(-\Delta)^s \zeta(x) dx -\int_{\R^n}  \left((-\Delta)^s\phi_s*\zeta(z)\right)\eta(z)  (-\Delta )^s  u(z)  dz\\
&\qquad+ \int_{\R^n}  \left((-\Delta)^s\phi_s*\zeta(z) \right)\eta(z)f(z)  dz\\
&=\int_{\R^n}u(x)(-\Delta)^s \zeta(x)dx -\int_{\R^n}\eta(z)\zeta(z) (-\Delta)^s u(z)  dz+ \int_{\R^n}\eta(z) f(z) \zeta(z) dz\\
&=\int_{\R^n}u(x)(-\Delta)^s \zeta(x)dx -\int_{\R^n}\zeta(z) (-\Delta)^s u(z)  dz+ \int_{\R^n} f(z) \zeta(z) dz,\\
\end{split}
\end{equation}
where in the third equation above we applied Fubini's theorem. The last equation follows since $\eta \equiv 1$ on the support of $\zeta$.  We thus infer
\begin{equation*}
\int_{\R^n} R_{\rho}(x) (-\Delta)^{s} \zeta(x) dx = \int_{\R^n} f(z) \zeta(z) dz,
\end{equation*}
which implies that $R_{\rho}(x)$ is a weak solution of
$$(-\Delta)^s R_{\rho}(x)= f\quad\mbox{in }B'_\rho.$$

To prove ii),
we first convolute (\ref{u_rep_local}) with $\gamma_{s-1,\epsilon}$ defined in Lemma \ref{lower_semi_cont} and then we apply Fubini's theorem to obtain
\begin{equation}\label{gamma_conv}
u*\gamma_{s-1,\epsilon}(x) =\int_{\R^n} \eta(\xi) \left(\phi_{s}*\gamma_{s-1,\epsilon}(x-\xi)\right) d\mu(\xi) +R_{\rho}(x)*\gamma_{s-1,\epsilon}.
\end{equation}
From Lemma \ref{lower_semi_cont} we know that $\gamma_{s-1,\epsilon}$ is a $C^{1,1}$ function. Therefore, taking the Laplacian in the $x$-variable we obtain
\begin{equation}\label{est3}
\Delta u*\gamma_{s-1,\epsilon}-\Delta R_{\rho}(x)*\gamma_{s-1,\epsilon}=\int_{\R^n} \eta(\xi) \left(\Delta \phi_{s}*\gamma_{s-1,\epsilon}(x-\xi)\right) d\mu.
\end{equation}
Moreover, since $-\Delta \phi_s$ is $(s-1)$-superharmonic, by Lemma (\ref{lower_semi_cont}) we also know that $ -\Delta \phi_s* \gamma_{s-1,\epsilon} \nearrow -\Delta \phi_s$,  and $\Delta u*\gamma_{s-1,\epsilon}-\Delta R_{\rho}(x)*\gamma_{s-1,\epsilon}$ converges to $\Delta u-\Delta R_{\rho}$ almost everywhere as $\epsilon \to 0^+$. Taking $\epsilon \to 0^+$ in (\ref{est3}) and applying the Monotone Convergence Theorem finishes the proof of ii).

To obtain an estimate for $R_{\rho}(x)$, we need to derive an explicit formula for it. First, we note that $\phi_{s}$ is the trace of $E_{s}$ on the thin space $\R^n$.
Also, we extend the cut-off function $\eta$ to the $(n+1)$-dimensional ball so that $\eta(x,y) \in C_0^{\infty}(B_{2\rho})$,  $\eta(x,y) \equiv 1$ on $B_{\rho}$ and $|\nabla \eta(x,y)| \leq \frac{1}{\rho}$ on $B_{2\rho}$. Given the extension result of Theorem \ref{extension*}, we know that
\begin{equation}\label{u_split}
\begin{split}
\phi_{s}*\mu_{\rho}(x)&=\int_{\R^n} \eta(\xi) \phi_{s}(x-\xi)
\left( (-\Delta)^s u(\xi) -f(\xi)\right)d\xi\\
&=\int_{\R^n} \eta (\xi) \phi_{s}(x-\xi) \left(\lim_{y \to 0^+} y^b \frac{\partial}{\partial y} \Delta_b U(\xi,y) \right)d\xi-\int_{\R^n}
\eta(\xi) \phi_s(x-\xi) f(\xi) d\xi\\
&=\int_{\R_+^{n+1}} y^b \Delta_b \left(\eta(\xi,y) E_{s}(x-\xi,y)\right) \Delta_b U(\xi,y) d\xi dy -\int_{\R^n}
\eta(\xi) \phi_s(x-\xi) f(\xi) d\xi\\
&=\int_{\R_+^{n+1}} y^b \Delta_b \left( E_{s}(x-\xi,y)\right) \Delta_b (\eta(\xi,y)U(\xi,y)) d\xi dy\\
& -\int_{\R_+^{n+1}} y^b
\left(U(\xi,y) \Delta_b \eta(\xi,y)+2\nabla U(\xi,y)\cdot \nabla \eta(\xi,y)\right) \Delta_b E_{s}(x-\xi,y) d\xi dy\\
&-\int_{\R_+^{n+1}} y^b \nabla \bigg(E_{s}(x-\xi,y) \Delta_b \eta(\xi,y) + 2 \nabla E_{s}(x-\xi,y) \cdot \nabla \eta(\xi,y) \bigg) \cdot \nabla U(\xi,y) d\xi dy\\
&-\int_{\R^n}
\eta(\xi) \phi_s(x-\xi) f(\xi) d\xi.\\
\end{split}
\end{equation}
The last equation above follows from  integration by parts. Combining (\ref{u_split}) with the definition  of $R_{\rho}$ in (\ref{R_rho_def}) we obtain

\begin{equation}\label{R}
R_{\rho}(x)=u(x)-\eta(x) u(x)+J_{\rho}(x),
\end{equation}
where
\begin{equation}\label{J}
\begin{split}
J_{\rho}(x):=&\int_{\R_+^{n+1}} y^b \left(U(\xi,y) \Delta_b \eta(\xi,y)+2\nabla U(\xi,y)\cdot \nabla \eta(\xi,y)\right) \Delta_b E_s(x-\xi,y) d\xi dy\\
&+\int_{\R_+^{n+1}} y^b \nabla \left(E_s(x-\xi,y) \Delta_b \eta(\xi,y) + 2 \nabla E_s(x-\xi,y) \cdot \nabla \eta(\xi,y) \right) \cdot \nabla U(\xi,y) d\xi dy\\
&+\int_{\R^n}\eta(\xi) \phi_s(x-\xi) f(\xi) d\xi.\\
\end{split}
\end{equation}
Clearly, $R_{\rho}=J_{\rho}$ on $B'_{\rho}$, and  $R_{\rho}(x)$ satisfies the estimate (\ref{R_bnd}).
\end{pf}
\section{$C^{1,1}_{\rm loc}$-regularity of  Problem (\ref{prb2})}\label{C1alpha_regularity}

We define
\begin{equation}\label{w}
w_{\rho}(x):=\phi_{s}*\mu_{\rho}(x)+ R_{\rho}(x).
\end{equation}
Recalling Lemma \ref{rep_local*}, we explicitly observe that
\begin{equation}\label{w_lap}
\Delta w_{\rho}(x):=\Delta \phi_{s}*\mu_{\rho}(x)+ \Delta R_{\rho}(x),
\end{equation}
and
$w_{\rho}=u$, \     $\Delta w_{\rho}=\Delta u$ for a.e.  $x\in B'_{\rho}$. From potential theory we know that $w_{\rho}$ and $-\Delta w_{\rho} $ are lower semi-continuous on $B'_{\rho}$. We begin with studying the behavior of $w_{\rho}$ on the support of the measure $\mu$.

\begin{lemma}\label{w_on_supp*}
For $x \in B'_{\rho} \cap \supp(\mu)$ we have $w_{\rho}(x)=\psi(x)$, where $\psi$ is the obstacle as in \eqref{prb2}.
\end{lemma}
\begin{pf}
Let $x_0 \in B'_{\rho}$ be such that $w(x_0)-\psi(x_0)=d>0$. We want to show that $x_0 \notin \supp(\mu)$. By the lower semi-continuity of $w$, there exists $r>0$ such that $w(x)-\psi(x)>{d}/{2}$ for all $x \in B'_r(x_0)$. Now, for any non-negative function $\phi \in C_0^{\infty}(B'_r(x_0))$  and strictly positive in $B'_{{r}/{2}}(x_0)$ such that $\|\phi\|_{L^{\infty}(B'_r(x_0))}<{d}/{2}$,  we see that the function $u-\phi$ is in the admissible set $\mathcal{B}$ defined in (\ref{prb2_admissible}).
Therefore,
$$\langle\mu,\phi\rangle=-\langle\mu, (u-\phi)-u\rangle<0.$$
This implies that $\mu(B'_{{r}/{2}}(x_0))=0$ which shows that $x_0 \notin \supp(\mu)$.
\end{pf}

\noindent
As a consequence, we obtain the following regularity result of $w_\rho$.

\begin{lemma}\label{w_reg*}
The function $w_{\rho}$ is in $C(\mathbb{R}^n)$, and it satisfies
\begin{equation}\label{w_bnd}
\|w_{\rho}\|_{L^{\infty}(B'_{\rho/2})}
\leq  C(\rho,n,s) \left(\|\psi\|_{L^{\infty}(\R^n)}+\|y^{{b}/{2}} U\|_{L^2(B_{2 \rho})}  +\| y^{{b}/{2}} \nabla U\|_{L^2(B_{2 \rho})}+\|f\|_{L^{\infty}(B'_{2 \rho})} \right).
\end{equation}
\end{lemma}

\begin{pf}
From Lemma \ref{w_on_supp*} we know that $w_{\rho}$, and hence $\phi_s* \mu_{\rho}$, is continuous on  $\supp(\mu_{\rho/2})$. Moreover, applying a heredity principle from potential theory (see \cite[p. 229]{SW}), we infer that $\phi_s* \mu_{\rho/2}$ is continuous on  $\supp(\mu_{\rho/2})$.
In addition,
\begin{equation}\label{w_bnd_part1}
\|\phi_s* \mu_{\rho/2}\|_{L^{\infty}(\supp(\mu_{\rho/2}))}\leq \|\phi_s* \mu_{\rho}\|_{L^{\infty}(\supp(\mu_{\rho/2}))}
\leq  C(\rho,n,s) \left(\|\psi\|_{L^{\infty}(\R^n)}+ ||R_{\rho}||_{L^{\infty}(B'_{\rho}})\right),
\end{equation}
where the first inequality follows from the heredity principle, and in the last inequality we used (\ref{u_rep_local}). Next, applying a continuity principle and a maximum principle from potential theory (see for instance \cite[p. 365]{L}), we conclude that $\phi_s* \mu_{\rho/2}$ is continous on $\mathbb{R}^n$, and
\begin{equation}\label{w_bnd_part2}
\|\phi_s* \mu_{\rho/2}\|_{L^{\infty}(\mathbb{R}^n)}
\leq  C(\rho,n,s) \left(\|\psi\|_{L^{\infty}(\R^n)}+ ||R_{\rho}||_{L^{\infty}(B'_{\rho})}\right).
\end{equation}
Finally, recalling (\ref{R_bnd}), we reach the desired result.
\end{pf}
\\
The following lower bound on $-\Delta w_{\rho}$ follows immediately from (\ref{w_lap}) and the positivity of the measure $\mu$.
\begin{corollary}\label{lapxw_upper_bnd*}
For every $B'_{\rho} \subset \mathbb{R}^n$,
\begin{equation}\label{lapxw_upper_bnd}
-\Delta w_{\rho} \geq -\|R_{\rho}
\|_{L^{\infty}(B_{\rho})}\quad\text {on } B'_{\rho}.
\end{equation}
\end{corollary}
We now obtain an upper bound  for $-\Delta w_{\rho}$ on the support of the measure $\mu$.
\begin{lemma}\label{lapxw_lower_bnd_support*}

 For $x_0 \in \supp (\mu_{\rho/2})$, we have
 \begin{equation}\label{lapxw_lower_bnd_support}
 -\Delta w_{\rho} (x_0) \leq -\Delta \psi (x_0),
 \end{equation}

\end{lemma}

\begin{pf}
In this proof, and for the sake of simplicity, we write $w$ for $w_{\rho}$ omitting the subscript $\rho$.
Let $x_0 \in \supp(\mu_{\rho/2})$.
From Lemma \ref{w_on_supp*} we know that $w(x_0)=\psi(x_0)$. Moreover, since $w$ is continuous and $w(x)=u(x)$ a.e. on $B'_{\rho}$, then $w(x) \geq \psi(x)$ on $B'_{\rho}$. Now we recall the definition of the usual mollifiers on $\R^n$

\begin{equation}\label{mollifier}
\omega_{\epsilon}(x):=\frac{1}{\epsilon^n} \omega\left(\frac{x}{\epsilon}\right),
\end{equation}
where
\begin{equation*} \omega(x)=\left\{\begin{array}{ll}
Ce^{\frac{1}{{|x}^2-1}} & |x|<1\\
0 & |x| \geq 1.\\
\end{array}\right.
\end{equation*}
We let $w_{\epsilon}(x)=\left(w(\cdot,0)*\omega_{\epsilon}\right)(x)$,  $\Delta w_{\epsilon}(x)=\left(\Delta_x w(\cdot,0)*\omega_{\epsilon}\right)(x)$, and  $\psi_{\epsilon}(x)=\psi*\omega_{\epsilon}(x)$.
Clearly, we also have $w_{\epsilon}(x)-\psi_{\epsilon}(x)\geq 0$ on $B'_{\rho}$.
For every $0<r<\rho-|x_0|$, applying Green's formula to $w_{\epsilon}-\psi_{\epsilon}$ on $B'_r(x_0)$ we obtain

\begin{equation*}
\begin{split}
w_{\epsilon}(x_0)-\psi_{\epsilon}(x_0)&=-\int_{B'_r(x_0)} \Delta (w_\epsilon-\psi_\epsilon)(x) G(x_0,x)dx +\frac{1}{|\partial{B'_r}(x_0)|}\int_{\partial B'_r(x_0)} (w_\epsilon-\psi_\epsilon)(x) dS\\
&\geq -\int_{B'_r(x_0)} \Delta (w_\epsilon-\psi_\epsilon)(x) G(x_0,x)dx.\\
\end{split}
\end{equation*}
From Corollary \ref{lapxw_upper_bnd*}, we know that $-(\Delta w (x)-\Delta \psi)$ is bounded from below on $B'_{\rho}$ by $a_{\rho}:=-\|R_{\rho}\|_{L^{\infty}(B_{\rho})}-\|\psi\|_{L^{\infty}(B_{\rho})}$.
This immediately implies that $-\Delta  (w_\epsilon-\psi_\epsilon)>a_{\rho}$ on $B'_{r}(x_0)$ for $0<r<\rho-|x_0|$.
Letting $\epsilon \to 0^+$ and applying Fatou's lemma we obtain
\begin{equation*}
-\int_{B'_r(x_0)}  \left(\Delta w(x)-\Delta \psi(x)\right) G(x_0,x)dx \leq w(x_0)-\psi(x_0)=0.
\end{equation*}
Thus, for every $0<r<\rho-|x_0|$, there exists a point $x_r \in B'_r(x_0)$ such that $-\Delta_x w(x_r)\leq -\Delta \psi(x_r)$. Choosing a sequence of points $x_r \to x_0$, by the lower semi-continuity of $-\Delta_x w$ we have
\begin{equation*}
-\Delta_x w(x_0)\leq \liminf_{r \to 0} -\Delta_x w(x_r) \leq \liminf_{r \to 0} -\Delta \psi(x_r) =-\psi(x_0)
\end{equation*}
and this finishes the proof of (\ref{lapxw_lower_bnd_support}).
\end{pf}
\\
With the results of Corollary \ref{lapxw_upper_bnd*} and  Lemma  \ref{lapxw_lower_bnd_support*}  in hand, adapting a maximum principle from potential theory, and then using the representation (\ref{w_lap}),  we can reach the following local bound on $\Delta w_{\rho}$ in $\mathbb{R}^n$.
\begin{lemma}\label{lapw_bnd*}
$\Delta w_{\rho}$ is locally bounded in $\mathbb{R}^n$.
 In particular,
\begin{equation}\label{lapw_bnd}
\begin{split}
 \|\Delta w_{\rho}\|_{L^{\infty}(B^+_{{\rho}/{2}})}  \leq  C(\rho,n,s) \bigg(&\|\Delta \psi\|_{L^{\infty}(\R^n)}+\|y^{{b}/{2}} U\|_{L^2(B_{2 \rho})}\\
 &+\left.\| y^{{b}/{2}} \nabla U\|_{L^2(B_{2 \rho})}+\|f\|_{L^{\infty}(B'_{2 \rho})} \right).
 \end{split}
\end{equation}
\end{lemma}
\begin{pf}
Combining the results of Corollary \ref{lapxw_upper_bnd*} and Lemma \ref{lapxw_lower_bnd_support*} we obtain the following estimate

\begin{equation}\label{lapw_bnd_support}
\begin{split}
\|\phi_s * \mu_{\rho}\|_{L^{\infty}(supp(\mu_{\rho}/2))}   \leq  C(\rho,n,s) \bigg(&\|\Delta \psi\|_{L^{\infty}(\R^n)}+\|y^{{b}/{2}} U\|_{L^2(B_{2 \rho})} \\ &+\left.\| y^{{b}/{2}} \nabla U\|_{L^2(B_{2 \rho})}+\|f\|_{L^{\infty}(B'_{2 \rho})} \right).
 \end{split}
\end{equation}
Applying the heredity principle and adapting a maximum principle from potential theory as we did in the proof of Lemma \ref{w_reg*}, we reach the estimate (\ref{lapw_bnd}).
\end{pf}
\\
The boundedness of $\Delta w_{\rho}$ immediately implies the following result.

\begin{corollary}\label{zero_measure}
For $x \in \R^n$, we have $\mu(\{x\})=0$.
\end{corollary}

\begin{pf}
Using (\ref{lapw_bnd}) in the representation (\ref {w_lap}) we see that
$\Delta \phi_{s}*\mu_{\rho}(x)$ is bounded in $B'_{\rho/2}$. Thus, for every $r<\rho/2-|x|$ we have
\begin{equation*}
\frac{1}{r^{n-2(s-1)}} \mu_{\rho}(B'_r)\leq C(s) \int_{B'_r} \Delta \phi_s (x-\xi) d\mu_{\rho}(\xi)<\infty.
\end{equation*}
Therefore,
$$\mu_{\rho}(B'_r) \leq C(s) r^{n-2(s-1)}.$$
Letting $r \to 0^+$, the result follows.
\end{pf}
\\
The following observation, which can be checked by direct computations, is crucial to obtain local bounds on all second order partial derivatives of $u$ on $B'_{\rho/2}$.

\begin{remark}\label{key_ineq*}
 $\phi_{s}$ satisfies
\begin{equation}\label{key_ineq}
\left|\frac{\partial}{\partial x_i x_j} \phi_{s} \right|\leq -C(n,s) \ \Delta \phi_{s}
\end{equation}
for all $i,j=1,\dots,n$.
\end{remark}
We can now prove the local boundedness of all second order partial derivatives of $u$ on $B'_{{\rho}/{2}}$. We first obtain the following representation:

\begin{lemma}\label{mixed_rep_local_upper*}
For all $i, j=1,\dots,n$ we have
\begin{equation}\label{mixed_rep_local_upper}
u_{x_i x_j} (x)= \phi_{x_i x_j} *\mu_{\rho}(x) + R_{x_i x_j} (x)
\end{equation}
for a.e. $x \in B'_{\rho}$.
\end{lemma}

\begin{pf}
For $i,j=1,\dots,n$, convoluting (\ref{u_rep_local}) in the tangential direction with $\gamma_{s-1,\epsilon}$, using Fubini's theorem, and then applying the operator $\frac{\partial}{\partial x_i \partial x_j}$, we get
\begin{equation}\label{gamma_extension_conv}
-u_{ x_i  x_j}*\gamma_{s-1,\epsilon}=\left(-(\phi_s)_{ x_i  x_j}*\gamma_{s-1,\epsilon}\right)*\mu_{\rho}-R_{ x_i  x_j}*\gamma_{s-1,\epsilon}.
\end{equation}
We have $(\phi_{s})_{x_i x_j} *\gamma_{s-1,\epsilon} \to (\phi_s)_{x_i x_j}$ pointwise on  $\R^n \setminus \{0\}$. By Corollary \ref{zero_measure}, we know that $\mu(\{0\})=0$. Therefore, the convergence  $(\phi_{s})_{x_i x_j}*\gamma_{s-1,\epsilon} \to (\phi_s)_{x_i x_j}$ is pointwise $\mu_{\rho}$-a.e. as $\epsilon \to 0^+$. Moreover, by (\ref{key_ineq}) and Lemma \ref{lower_semi_cont} we see that
\begin{equation*}
|(\phi_{s})_{x_i x_j} *\gamma_{s-1,\epsilon}|\leq -C\Delta \phi*\gamma_{s-1,\epsilon}\nearrow -C \Delta \phi_s.
\end{equation*}
From Lemma \ref{lapw_bnd*} and the representation (\ref {w_lap}), we know that $-C \Delta \phi_s * \mu_{\rho}<\infty$. Therefore, an application of the Dominated Convergence Theorem yields (\ref{mixed_rep_local_upper*}).
\end{pf}
\\
The following corollary follows immediately from Lemma \ref{lapw_bnd*}, Remark \ref{key_ineq*}, and Lemma \ref{mixed_rep_local_upper*}.

\begin{corollary}\label{C11_reg*}
$u$ is in $C^{1,1}_{\rm loc}(\mathbb{R}^n)$. In particular,
\begin{equation} \label{C1alpha_reg}
\begin{split}
\|u\|_{C^{1,1}( B'_{{\rho}/{2}})}\leq  C(n,s,\rho) \bigg(&\|\Delta \psi\|_{L^{\infty}(\R^n)}+\|y^{b/2} U\|_{L^2(B_{2 \rho})}\\  &+\left.\| y^{b/2} \nabla U\|_{L^2(B_{2 \rho})}+\|f\|_{L^{\infty}(B'_{2 \rho})} \right).
\end{split}
\end{equation}
\end{corollary}

\section{$H^{1+s}$-regularity of Problem (\ref{prb2})}\label{H1+s_regularity}
Following  the approach in \cite[Lemma 1.9]{Sc1}, next we  show that the solution to Problem (\ref{prb2}) is in $  H^{1+s}_{\rm loc}(\R^n)$.

\begin{theorem}\label{H1+s_reg*}
Let $u$ be the minimizer of $I[\cdot]$. Then, $u \in H_{\rm loc}^{1+s}(\R^n)$.  In particular,
\begin{equation}\label{H1+s_reg}
\|u\|_{H^{1+s}(B'_{{\rho}/{4}})}\leq C(\rho) \left(\|\Delta \psi\|_{L^{\infty}(\R^n)}+\|y^{{b}/{2}} U\|_{L^{2}(B_{2\rho})}+\|y^{{b}/{2}} \nabla U\|_{L^{2}(B_{2\rho})} + \|f\|_{L^{\infty}(B'_{2\rho})}\right).
\end{equation}
\end{theorem}

\begin{pf} We define
\begin{equation}\label{Q_beta}
Q_{\beta}(x):=\frac{1}{|x|^{n-\beta}},\quad x \in \R^n,
\end{equation}
which is equal (up to a constant) to the Riesz kernel of order $\beta$. It is a classical fact from potential theory (see, for example, \cite[Theorem 1.20]{L}) that the following relation between the $L^2$-norm of the Riesz potential of order $\beta/2$ of a compactly supported measure $\nu$ and the Riesz energy of order $\beta$ of the same measure holds:
\begin{equation}\label{measure_energy}
\int_{\R^n} | Q_{\beta/2}* \nu|^2 dx =C(n,\beta)\int_{\R^n} (Q_{\beta}* \nu) d\nu\quad\text{for all } 0<\beta<{n}/{2}.
\end{equation}
Setting $\nu=\mu_{\rho/4}$ and $\beta=2(s-1)$ in (\ref{measure_energy}) we obtain
\begin{equation}\label{formula2}
\begin{split}
\|Q_{s-1}*&\mu_{\rho/4}\|_{L^2(\mathbb{R}^n)}=C(n,s) \int_{\R^n} Q_{2(s-1)}*\mu_{\rho/4} d\mu_{\rho/4}\\
&=-C(n,s)\int_{\R^n} \Delta \phi_s *\mu_{\rho/4} d\mu_{\rho/4}\\
&\leq C(n,s,\rho) \| \Delta \phi_s *\mu_{\rho/4} \|_{L^{\infty}(B'_{\rho/2})}\\
&\leq C(n,s,\rho) \| \Delta \phi_s *\mu_{\rho/2} \|_{L^{\infty}(B'_{\rho/2})}\\
&\leq C(n,s,\rho) \left(\|\Delta \psi\|_{L^{\infty}(\R^n)}+\|y^{{b}/{2}} U\|_{L^{2}(B_{2\rho})}+\|y^{{b}/{2}} \nabla U\|_{L^{2}(B_{2\rho})}+ \|f\|_{L^{\infty}(B'_{2\rho})}\right),\\
\end{split}
\end{equation}
where the last inequality follows from Lemma \ref{lapw_bnd*}.
Noticing that
\begin{equation*}
\begin{split}
- \widehat{(-\Delta)^{\frac{s-1}{2}} \Delta \phi_s}&=C(n,s) \widehat{(-\Delta)^{\frac{s-1}{2}} Q_{2(s-1)}}\\
&=C(n,s) |\xi|^{s-1}\frac{1}{|\xi|^{2s-2}}=C(n,s) \frac{1}{|\xi|^{s-1}}=C(n,s) \widehat{Q_{s-1}},
\end{split}
\end{equation*}
we conclude that
\begin{equation*}
(-\Delta)^{\frac{s-1}{2}}\Delta \phi_s * \mu_{{\rho}/{4}} \in L^2(B'_{\rho/2}),
\end{equation*}
and
\begin{equation}\label{lapu_diff}
\begin{split}
\|(-\Delta)^{\frac{s-1}{2}}\Delta \phi_s * \mu_{{\rho}/{4}}\|_{L^2(B'_{\rho/2})}\leq C(\rho) &\Bigl(\|\Delta \psi\|_{L^{\infty}(\R^n)}\\
&+\|y^{{b}/{2}} U\|_{L^{2}(B_{2\rho})}+\|y^{{b}/{2}} \nabla U\|_{L^{2}(B_{2\rho})}+ \|f\|_{L^{\infty}(B'_{2\rho})}\Bigr).
\end{split}
\end{equation}
Combining (\ref{lapu_diff}) with the representation (\ref{lapxu_rep_local}), we infer that
$\Delta u \in H^{s-1}_{\rm loc}(\R^n)$ and
\begin{equation}\label{H1+s_reg}
\|\Delta u\|_{H^{s-1}(B'_{{\rho}/{4}})}\leq C(\rho) \left(\|\Delta \psi\|_{L^{\infty}(\R^n)}+\|y^{{b}/{2}} U\|_{L^{2}(B_{2\rho})}+\|y^{{b}/{2}} \nabla U\|_{L^{2}(B_{2\rho})} + \|f\|_{L^{\infty}(B'_{2\rho})}\right).
\end{equation}
Moreover, we notice that, for all $i,j:1,..,n$,
\begin{equation}\label{mixed_diff}
\int_{B'_{\rho/8}}  \left|\widehat{(-\Delta)^{\frac{s-1}{2}} (\phi_s)_{x_i x_j}}\right|^2=\int_{B'_{\rho/8}}   |\xi|^{2(s-1)} \left|\widehat{(\phi_s)_{x_i x_j}}\right|^2 \leq \int_{B'_{\rho/8}}  |\xi|^{2(s-1)} \left|\widehat{\Delta \phi_s} \right|^2 ,
\end{equation}
where the last inequality follows from Remark \ref{key_ineq*} and the fact that $\Delta \phi_s \leq 0$.
Combining (\ref{mixed_diff}) with  (\ref{lapu_diff}) and the representation (\ref{mixed_rep_local_upper}) with $y=0$, we reach the desired conclusion.
\end{pf}

It is easy to see that the arguments of Section \ref{C1alpha_regularity} and Section \ref{H1+s_regularity} can be applied to the minimizer $u_0$ of $I_0[\cdot]$ as well. Therefore, we have the following regularity of $u_0$.

\begin{theorem}\label{u0_regularity}
The minimizer $u_0$ of $I_0[\cdot]$ is in $C_{\rm loc}^{1,1} (\R^n)\cap H^{1+s}_{\rm loc}(\R^n)$.
\end{theorem}

\section{Regularity of the extension problem (\ref{extension_prb})}
In this section, we extend the regularity results of Corollary
\ref{C11_reg*} and Theorem
\ref{H1+s_reg*} of problem (\ref{prb2}) to the extension problem (\ref{extension_prb}). We start with extending
the representations obtained in Lemma \ref{rep_local*} to the upper half space as follows. We recall that $E_s(x,y)$ is the fundamental solution of the b-biLaplace
equation as defined in \eqref{fundamental}.
\begin{lemma}\label{rep_local_upper*}
For every $y>0$, and  every ball $B'_{\rho}$ in  $\R^n$, the b-biharmonic extension $U$ of $u$ and its Laplacian admit the following representations:
\begin{enumerate}[i)]
\item For  every $y>0$ and almost every $x\in  B'_{\rho}$
\begin{equation}\label{u_rep_local_upper}
U(x,y)=\int_{\R^n} E_s(x-\xi,y) \ d\mu_{\rho}(\xi)+L_{\rho}(x,y),
\end{equation}
where $L_\rho(x,y)$ is the b-biharmonic extension of $R_\rho(x)$ to the upper half-space.
\item For  every $y>0$ and almost every $x\in  B'_{\rho}$
\begin{equation}\label{lapu_rep_local_upper}
\Delta U(x,y)=\int_{\R^n} \Delta E_s(x-\xi,y) \ d\mu_{\rho}(\xi)+\Delta L_{\rho}(x,y).
\end{equation}
\item For almost every $x\in  B'_{\rho}$
\begin{equation}\label{uyy_rep_local}
U_{yy}(x,0)= C(n,s,\rho) \Delta \phi _s *\mu_{\rho}(x) + (L_\rho)_{yy}(x,0).
\end{equation}
\item $L_{\rho}$ satisfies the estimate
\begin{equation}\label{L_bnd}
\|L_{\rho}\|_{W^{2,\infty}(B_{\rho})} \leq C(n,s,\rho) \left(\|y^{{b}/{2}} U\|_{L^2(B_{2 \rho})}  +\| y^{{b}/{2}} \nabla U\|_{L^2(B_{2 \rho})}+\|f\|_{L^{\infty}(B'_{2 \rho})} \right).
\end{equation}
\end{enumerate}
\end{lemma}
\begin{pf}
Convoluting (\ref{u_rep_local}) in the tangential directions with the poisson kernel $P_s(x,y)$ given in (\ref{Poisson_kernel}) we obtain $i)$. To prove $ii)$,  for every fixed $y>0$, we convolute (\ref{u_rep_local_upper}) in the tangential directions with $\gamma_{s-1,\epsilon}$ and then apply Fubini's theorem to get
\begin{equation*}\label{gamma_conv}
U(\cdot,y)*\gamma_{s-1,\epsilon}=\left(E_s(\cdot,y)*\gamma_{s-1,\epsilon}\right)*\mu_{\rho}+L_\rho(\cdot,y)*\gamma_{s-1,\epsilon},
\end{equation*}
for every $y>0$ and a.e. $x$ in $B'_{\rho}$. Applying the operator $-\Delta$ we obtain
\begin{equation}\label{gamma_extension_conv1}
-\Delta U(\cdot,y)*\gamma_{s-1,\epsilon}=\left(-\Delta E_s(\cdot,y)*\gamma_{s-1,\epsilon}\right)*\mu_{\rho}-\Delta  L_\rho (\cdot,y) *\gamma_{s-1,\epsilon},
\end{equation}
for every $y>0$ and a.e. $x$ in $B'_{\rho}$.
Noticing that $0\leq -\Delta E_s \leq -\Delta_x \phi_s$, and recalling that $-\Delta_x \phi_s$ is $(s-1)$-superharmonic, we infer that for every $y>0$ and a.e. $x \in B'_{\rho}$,
$$
\Delta E_s(\cdot,y)*\gamma_{s-1,\epsilon} \nearrow \Delta E_s (\cdot,y) \ \text { \ as } \epsilon \to 0^+.
$$
Letting $\epsilon \to 0^+$ in (\ref{gamma_extension_conv1}) and applying the Monotone Convergence Theorem finishes the proof of $ii)$. Finally, letting $y \to 0^+$ in (\ref{lapu_rep_local_upper}) and recalling (\ref{lapxu_rep_local}) yields $iii)$.

To prove $iv)$, we carry similar computations as in the proof of part $iii)$ of Lemma \ref{rep_local*}.
From (\ref{u_rep_local_upper}) we have
\begin{equation}\label{L_rho}
L_{\rho}(x,y)=U(x,y)-\int_{\R^n} E_s(x-\xi,y) \ d\mu_{\rho}(\xi),
\end{equation}
for every $y>0$ and almost every $x\in  B'_{\rho}$.
Now we compute
\begin{equation*}
\begin{split}
\int_{\R^n} E_s&(x-\xi,y) \ d\mu_{\rho}(\xi)=\int_{\R^n} \eta(\xi) E_{s}(x-\xi,y)
\left( (-\Delta)^s u(\xi) -f(\xi)\right)d\xi\\
&=\int_{\R^n} \eta (\xi) E_{s}(x-\xi,y) \left(\lim_{z \to 0^+} z^b \frac{\partial}{\partial z} \Delta_b U(\xi,z) \right)d\xi-\int_{\R^n}
\eta(\xi) E_s(x-\xi,y) f(\xi) d\xi\\
&=\int_{\R_+^{n+1}} z^b \Delta_b \left(\eta(\xi,z) E_{s}(x-\xi,y-z)\right) \Delta_b U(\xi,z) d\xi dz -\int_{\R^n}
\eta(\xi) E_s(x-\xi,y) f(\xi) d\xi\\
&=\int_{\R_+^{n+1}} z^b \Delta_b \left( E_{s}(x-\xi,y-z)\right) \Delta_b (\eta(\xi,z)U(\xi,z) d\xi dz\\
& -\int_{\R_+^{n+1}} z^b
\left(U(\xi,z) \Delta_b \eta(\xi,z)+2\nabla U(\xi,z)\cdot \nabla \eta(\xi,z)\right) \Delta_b E_{s}(x-\xi,y-z) d\xi dz\\
&-\int_{\R_+^{n+1}} z^b \nabla \bigg(E_{s}(x-\xi,y-z) \Delta_b \eta(\xi,z) + 2 \nabla E_{s}(x-\xi,z) \cdot \nabla \eta(\xi,z) \bigg) \cdot \nabla U(\xi,z) d\xi dz\\
&-\int_{\R^n}
\eta(\xi) E_s(x-\xi,y) f(\xi) d\xi.\\
\end{split}
\end{equation*}
Noticing that the first term on the right hand side above is equal to $\eta(x,y)U(x,y)$,
we see that for $(x,y) \in B_{\rho}$,
\begin{equation}\label{L_rho_local}
\begin{split}
L_{\rho}&(x,y)= \int_{\R_+^{n+1}} z^b
\left(U(\xi,z) \Delta_b \eta(\xi,z)+2\nabla U(\xi,z)\cdot \nabla \eta(\xi,z)\right) \Delta_b E_{s}(x-\xi,y-z) d\xi dz\\
&+\int_{\R_+^{n+1}} z^b \nabla \bigg(E_{s}(x-\xi,y-z) \Delta_b \eta(\xi,z) + 2 \nabla E_{s}(x-\xi,z) \cdot \nabla \eta(\xi,z) \bigg) \cdot \nabla U(\xi,z) d\xi dz\\
&+\int_{\R^n}
\eta(\xi) E_s(x-\xi,y) f(\xi) d\xi,\\
\end{split}
\end{equation}
which clearly satisfies the estimate
(\ref{L_bnd}).
\end{pf}

\noindent
Now we prove the following local $C^{1,1}$-regularity result for the extension problem (\ref{extension_prb}).
\begin{theorem}\label{C11_reg_upper*}
$U$ is in $C^{1,1}_{\rm loc}(\overline{\mathbb{R}^{n+1}})$. In particular,
\begin{equation} \label{C1alpha_reg}
\begin{split}
\|U\|_{C^{1,1}(B^+_{\rho/2} \cup \  B'_{{\rho}/{2}})}\leq  C(n,s,\rho) &\bigg(\|\Delta \psi\|_{L^{\infty}(\R^n)}+\|y^{b/2} U\|_{L^2(B_{2 \rho})}\\ &+\left.\| y^{b/2} \nabla U\|_{L^2(B_{2 \rho})}+\|f\|_{L^{\infty}(B'_{2 \rho})} \right).
\end{split}
\end{equation}

\end{theorem}
\begin{pf}
Letting $y \to 0$ in (\ref{u_rep_local_upper}) we infer that $U$ is continuous in $\overline{\mathbb{R}^{n+1}}$. Moreover, for $(x,y) \in B_{\rho/2}$, the representation  (\ref{u_rep_local_upper}) yields
\begin{equation*}
\begin{split} |U(x,y)|&\leq \int_{\R^n} E_s(x-\xi,y) \ d\mu_{\rho/2}(\xi)+\left|L_{\rho/2}(x,y)\right|\\
 &\leq \|\phi_{s}*\mu_{\rho/2}(x)\|_{L^{\infty}(B'_{\rho/2})} + \|L_{\rho/2}\|_{L^{\infty}(B_{\rho})}, \\
\end{split}
\end{equation*}
where the second inequality can be explicitly verified by comparing (\ref{fundamental_fractional}) and (\ref{fundamental}).
Employing the estimates (\ref{w_bnd}) and $(\ref{L_bnd})$, we reach
\begin{equation}\label{U_bnd}
\begin{split}
\|U\|_{L^{\infty}(B^+_{\rho/2} \cup B'_{\rho/2})}
\leq  C(\rho,n,s) &\bigg(\|\psi\|_{L^{\infty}(\R^n)}+\|y^{{b}/{2}} U\|_{L^2(B_{2 \rho})} \\ &+\left.\| y^{{b}/{2}} \nabla U\|_{L^2(B_{2 \rho})}+\|f\|_{L^{\infty}(B'_{2 \rho})} \right).
\end{split}
\end{equation}
Similarly, from the representations (\ref{lapu_rep_local_upper}) and \eqref{uyy_rep_local}, we infer that for $(x,y) \in B_{\rho/2}$
\begin{equation*}
\begin{split}
|\Delta U(x,y)|&\leq -\int_{\R^n} \Delta E_s(x-\xi,y) \ d\mu_{\rho/2}(\xi)+|\Delta L_{\rho/2}(x,y)|\\
&\leq \|\Delta \phi_{s}*\mu_{\rho/2}(x)\|_{L^{\infty}(B'_{\rho/2})}+|\Delta L_{\rho/2}(x,y)|,
\end{split}
\end{equation*}
and
\begin{equation*}
|U_{yy}(x,0)|\leq  C(n,s,\rho) \|\Delta \phi_{s}*\mu_{\rho/2}(x)\|_{L^{\infty}(B'_{\rho/4})} + |(L_{\rho/2})_{yy}(x,0)|.
\end{equation*}
Using Corollary \ref{C11_reg*} and applying the estimate (\ref{L_bnd}), we reach the desired result.
\end{pf}

Finally, from the regularity result of Theorem \ref{H1+s_reg*}, we  deduce the $H_{\rm loc}^{3}(\R^{n+1},y^b)$-regularity of the extension function $U$ of $u$ to the upper half-space. This regularity result will be crucial in forthcoming work on a time dependent thin obstacle problem for the b-biLaplace operator.

\begin{theorem}\label{gradLapU_reg*}
Let $u$ be the minimizer of $I[\cdot]$  and let $U$ be the extension of $u$ to the upper half space $\R^+_{n+1}$. Then, $\Delta_b U\in H_{\rm loc}^{1}(\R^{n+1},y^b)$. Moreover,
\begin{equation}\label{H1+s_to_gradLapU}
\begin{split}
\|\nabla (\Delta_b U)\|_{L^2(B_{\rho/8}, y^b)} \leq C(\rho) &\left(\|\Delta \psi\|_{L^{\infty}(\R^n)}+\|y^{{b}/{2}} U\|_{L^{2}(B_{2\rho})}+\|y^{{b}/{2}} \nabla U\|_{L^{2}(B_{2\rho})}\right.\\
&\left.
+ \|f\|_{L^{\infty}(B'_{2\rho})}\right).
\end{split}
\end{equation}
\end{theorem}

\begin{pf}
From Theorem \ref{H1+s_reg*}, we know that $\Delta u \in H^{s-1}_{\rm loc}(\R^n)$. Therefore, by the Caffarelli-Silvestre Extension Theorem, $\Delta u$ can be extended to the upper half-space by a function $W \in H^1_{\rm loc}(\R^{n+1},y^b)$ which is a solution of
\begin{equation*}
\left\{\begin{array}{ll}
\Delta_b W=0 & \text { in } \R^{n+1}_+\\
W_y=(-\Delta)^s u& \text { in } \R^{n}\\
\end{array}
\right.
\end{equation*}
 in the weak sense
\begin{equation}\label{W_extension}
\int_{B_{r}^+} y^b \ \nabla W \cdot \nabla \zeta =\langle \Delta u,\zeta\rangle_{\dot{H}^{\gamma}(\mathbb{R}^n)},
\end{equation}
for all $\zeta \in C_0^{\infty}(B_{r})$. To reach the first conclusion of the theorem, we show that $W=-\Delta_b U$ on $B_{r}$. This is equivalent to proving
\begin{equation}\label{to_show}
\int_{B_{r}} y^b (W+\Delta_b U) \psi=0 \quad\text{for all }\psi \in C_0^{\infty}(B_{r}).
\end{equation}
We observe that for such a function $\psi$  the problem
 \begin{equation}\label{problem}
\left\{
\begin{array}{ll}
\Delta_b \zeta=\psi &\text { in } B_{r}^+\\
\zeta=0 & \text { on } (\partial B_{r})^+\\
\zeta_y=0 & \text { on } B_{r}'\\
\end{array}
\right.
\end{equation}
admits a unique solution $\zeta_0$ whose even reflection to the whole ball $B_{r}$ is in  $C_0^{\infty}(B_{r})$. Substituting $\zeta_0$ for $\xi$ in (\ref{W_extension}) and integrating the left-hand-side by parts we obtain
\begin{equation}\label{identity_1}
\int_{B_{r}^+} y^b  W  \Delta_b \zeta_0 =-\langle \Delta u,\zeta_0\rangle_{\dot{H}^{\gamma}(\mathbb{R}^n)}.
\end{equation}

\noindent
Notice that, with the regularity result of Theorem \ref{H1+s_reg*} in hand, both the minimizer  $u$ and its extension $U$ satisfy
\begin{equation}\label{identity_2}
\int_{B_{r}^+} y^b  \Delta_b U  \Delta_b \zeta_0 =\langle \Delta u,\zeta_0\rangle_{\dot{H}^{\gamma}(\mathbb{R}^n)}.
\end{equation}
Combining (\ref{identity_1}) and (\ref{identity_2}) proves (\ref{to_show}).

For the second part of the theorem, set $r=\rho/4$, and let $\eta \in C_0^{\infty}(B_{\rho/4})$ such that $\eta \equiv 1$ on $B_{\rho/8}$.
Plugging in $\eta^2 \Delta_b U$ for $\zeta$ in (\ref{W_extension}) we obtain
\begin{equation*}
\begin{split}
\int_{B_{\rho/4}} y^b \ \nabla (\Delta_b U) \cdot \nabla (\eta^2 \Delta_b U)= \langle \Delta u,\eta^2 \Delta u\rangle_{\dot{H}^{\gamma}(\mathbb{R}^n)}
&= \int_{\R^n} |\xi|^{2(s-1)} \widehat{\Delta u} \ \widehat{\eta^2 \Delta u}\\
&\leq C(n,\rho) \|\Delta u\|^2_{H^{s-1}(B'_{\rho/4})}.\\
\end{split}
\end{equation*}
Expanding the left hand side we get
\begin{equation*}
\int_{B_{\rho/4}^+} y^b \ \eta^2 | \nabla (\Delta_b U)|^2 \leq C(n,\rho) \|\Delta u\|^2_{H^{s-1}(B'_{\rho/4})}  -2 \int_{B_{\rho/4}} y^b \eta (\Delta_b U) \nabla \eta \cdot \nabla (\Delta_b U).
\end{equation*}
Applying Young's inequality we infer
\begin{equation*}
\|\nabla (\Delta_b U)\|_{L^2(B_{\rho/8}, y^b)} \leq C(n,\rho) \left(\|\Delta u\|^2_{H^{s-1}(B'_{\rho/4})} +\|\Delta_b U\|_{L^2(B_{\rho/4}, y^b)} \right).
\end{equation*}
Finally, employing  Theorem  \ref{H1+s_reg*}, we reach (\ref{H1+s_to_gradLapU}).
\end{pf}

\begin{corollary}\label{H3_reg*}
Let $u$ be the minimizer of $I[\cdot]$  and let $U$ be the extension
of $u$ to the upper half space $\R_+^{n+1}$. Then, $U\in
H_{\rm loc}^{3}(\R^{n+1},y^b)$. Moreover,
\begin{equation}\label{H1+s_to_H3}
\begin{split}
\|U\|_{H^3(B_{\rho/8}, y^b)}\leq C(\rho) \bigg(&\|\Delta \psi\|_{L^{\infty}(\R^n)}+\|y^{{b}/{2}} U\|_{L^{2}(B_{2\rho})}+\|y^{{b}/{2}} \nabla U\|_{L^{2}(B_{2\rho})}\\
&\left.+\|f\|_{L^{\infty}(B'_{2\rho})}\right),
\end{split}
\end{equation}
for all $0<\rho<1$.
\end{corollary}

\begin{pf}
From Theorem \ref{gradLapU_reg*} and the representation (\ref{lapu_rep_local_upper}), we infer
$$\int_{\R^n} \nabla \left(\Delta E_s(x-\xi,y)\right) d\mu_{\rho/8}(\xi)\in L^2(B_{{\rho}/{8}},y^b),
$$
with the corresponding norm controlled by the right-hand-side of (\ref{H1+s_to_H3}).
 Moreover, we can easily check by direct computations that
 $$
 |(E_s)_{x_i x_j x_k}|\leq C(n,s) |\nabla (\Delta E_s)|$$
  for all $i,j,k=1,\dots,n+1$. Therefore,
 \begin{equation}\label{H1+s_to_H3_est}
 \int_{\R^n} (E_s)_{x_i x_j x_k}(x-\xi,y)d\mu_{\rho/2}(\xi) \in L^2(B_{{\rho}/{8}},y^b),
 \end{equation}
 with its norm controlled by the right-hand-side of (\ref{H1+s_to_H3}).
 Using this result in  the representation (\ref{u_rep_local_upper}), we reach the desired conclusion.
\end{pf}

{\bf Data Availability:}
Data sharing not applicable to this article as no datasets were generated or analysed during the current study.

\end{document}